\documentclass[12pt]{article} \usepackage{verbatim}
\usepackage{graphicx}
\usepackage{hyperref,url}
\usepackage{amsmath,amsthm,amssymb}
\newcommand{\nocopyright}{
No Copyright\thanks{
The authors hereby waive all copyright
and related or neighboring rights to this work,
and dedicate it to the public domain.
This applies worldwide.
}}
\title{Equivariant division}
\author{Prajeet Bajpai \and Peter G. Doyle}
\date{Version dated 13 April 2016
\\ \nocopyright
}

\newtheorem{theorem}{Theorem}

\newtheorem{prop}[theorem]{Proposition}

\newtheorem{problem}{Problem}

\newcommand{\union}{\cup}
\newcommand{\cross}{\times}

\newcommand{\id}{{\mbox{id}}}

\newcommand{\then}{\lhd}

\newcommand{\fabc}{f_{\alpha,\beta,\gamma}}
\newcommand{\hab}{h_{\alpha,\beta}}

\newcommand{\exstart}{\begin{gathered}}
\newcommand{\exend}{\end{gathered}}
\newcommand{\row}[1]{|_{#1}}
\newcommand{\FP}{\mathrm{FP}}
\newcommand{\bp}{*}
\begin{document}

\maketitle

\begin{abstract}

Let $C$ be a non-empty finite set, and 
$\Gamma$ a subgroup of the symmetric group $S(C)$.
Given a bijection $f:A \cross C \to B \cross C$,
the problem of \emph{$\Gamma$-equivariant division} is to 
find a \emph{quotient bijection} $h:A \to B$
respecting whatever symmetries $f$ may have under 
the action of $S(A) \cross S(B) \cross \Gamma$.
Say that $\Gamma$ is \emph{fully cancelling} if this is possible
for any $f$,
and \emph{finitely cancelling} if it is possible providing $A,B$ are finite.
Feldman and Propp showed that
a permutation group is finitely cancelling
just if it has a globally fixed point.
We show here that
a permutation group is fully cancelling
just if it is trivial.
This sheds light on the fact that all known division algorithms that eschew
the Axiom of Choice depend on fixing an ordering for the elements of $C$.
\end{abstract}

\section{Introduction}
Let $C$ be a non-empty finite set, and 
$\Gamma$ a subgroup of the symmetric group $S(C)$.
Given a bijection $f:A \cross C \to B \cross C$,
the problem of \emph{$\Gamma$-equivariant division} is to 
find a \emph{quotient bijection} $h:A \to B$
respecting whatever symmetries $f$ may have under 
the action of $S(A) \cross S(B) \cross \Gamma$.

Specifically,
given
\[
(\alpha,\beta,\gamma) \in S(A) \cross S(B) \cross \Gamma
,
\]
let
\[
\fabc = (\alpha^{-1} \cross \gamma^{-1}) \then f \then (\beta \cross \gamma)
,
\]
and
\[
\hab = \alpha^{-1} \then h \then \beta
,
\]
where the symbol $\then$, pronounced `then',
represents the composition of functions
in the natural order, with first things first:
\[
(p \then q)(x)=q(p(x))
.
\]
We say that $h$ is a \emph{$\Gamma$-equivariant quotient} of $f$
if whenever $\fabc = f$ we have $\hab = h$.
$\Gamma$ is \emph{fully cancelling} if
every bijection $f:A \cross C \to B \cross C$
has a $\Gamma$-equivariant quotient,
and \emph{finitely cancelling} if this is true providing $A,B$ are finite.

Feldman and Propp
\cite{fp}
looked at the finite case.
They showed that
the subgroup $S(C,\star)$ fixing a designated basepoint $\star \in C$
is finitely cancelling,
but unless $C$ is a singleton, the full group $S(C)$ is
not.
Going further, they gave a beautiful proof that
$\Gamma$ is finitely cancelling just if it has a globally fixed point.

Here we are interested in the infinite case.
The general problem of division is to produce
from $f:A \cross C \to B \cross C$
\emph{any} quotient bijection $h:A \to B$,
equivariant or not.
Known division methods
that eschew the Axiom of Choice
(cf. \cite{conwaydoyle:three,doyleqiu:four,schwartz:four})
produce quotients
that respect any symmetries under the action of $S(A) \cross S(B)$,
so they are at least $S_0(C)$-equivariant, where $S_0(C)$ is the
trivial subgroup of $S(C)$.
But these methods depend on fixing an ordering of $C$, suggesting that
this is the most equivariance we can hope for.
And indeed, we will show that $\Gamma$ is fully cancelling just if it
is the trivial subgroup $S_0(C)$.

\section{Finitely cancelling}

For starters,
Feldman and Propp showed that if you specify a base point $\bp \in C$,
the subgroup $S(C,\bp)$ of $S(C)$ that fixes $\bp$
is finitely cancelling.

Here's the argument.
For $c \in C$ define a map (not generally a bijection)
\[
f\row{c}:A \to B
,
\]
\[
f\row{c}(a) = f \then \pi_1(a,c)
,
\]
where
\[
\pi_1((x,y)) = x
.
\]
Let
\[
p(a) = f\row{\bp}(a) = (f \then \pi_1)(a,\bp)
\]
and
\[
q(b) = f^{-1}\row{\bp}(b) = (f^{-1} \then \pi_1)(b,\bp)
.
\]
Because $A$ is finite, the composition $p \then q$ has some cycles.
Let $X \subset A$ be the union of all these cycles.
The restriction $p|X$ is a partial bijection from $A$ to $B$.
Subtract $p|X \cross \id_C$ from $f$
(cf. \cite{conwaydoyle:three,doyleqiu:four,fp})
to get a bijection from $(A-X) \cross C$ to
$(B-p(X)) \cross C$.
Proceed by recursion to get a bijection $\FP(f,\bp): A \to B$.

To sum up:

\begin{prop}[Feldman-Propp] \label{star}
If some $\bp \in C$ is fixed by every
$g \in \Gamma$, $\Gamma$ is finitely cancelling.
\end{prop}

We can collect the various bijection $\FP(f,c)$ for $c \in C$
into a new bijection
\[
\bar{f}:A \cross C \to B \cross C
,
\]
\[
\bar{f}((a,c)) = (\FP(f,c)(a),c)
.
\]
This new bijection $\bar{f}$ satisfies
\[
\bar{f}((a,c)) = (\bar{f}\row{c}(a),c)
.
\]
We will call any bijection that preserves the second coordinate in this
way a \emph{parallel bijection}.

By combining all the bijections $\FP(f,c)$ in this parallelization $\bar{f}$,
we obviate the need to choose a basepoint, so Proposition \ref{star}
implies (and follows from):
\begin{prop} \label{parallel}
To a finite bijection
$f: A \cross C \to B \cross C$
we can associate in a fully equivariant manner a new bijection $\bar{f}$ with
\[
\bar{f}((a,c)) = (\bar{f}_c(a),c)
\]
where $\bar{f}_c:A \to B$ is a bijection for each $c \in C$.
\end{prop}

In light of Proposition \ref{parallel},
$\Gamma \subset S_C$ is finitely cancelling just if
any finite parallel bijection has a $\Gamma$-equivariant quotient.
Indeed, to any finite $f$ we can associate its parallelization $\bar{f}$;
if $\bar{f}$ has a $\Gamma$-equivariant quotient then so does $f$;
if it does not, then $\Gamma$ is not cancelling.

This does not necessarily mean that
in every finite division problem we can safely parallelize
$f$ as our first step.
It could be that $f$ has a $\Gamma$-equivariant quotient
while its parallelization $\bar{f}$
does not.
(See \ref{prob:parallel}.)

Proposition \ref{parallel} fails in the infinite case;
this fact underlies the counterexamples we will produce there.

\section{Not finitely cancelling}

We begin with counterexamples in the finite case,
all obtained using the method of Feldman and Propp.

The simplest case is $C=\{a,b\}$.
Take $A = \{x,y\}$, $B=\{1,2\}$, and
\[
\exstart
f =
\begin{array}{l|ll}
&x&y\\
\hline
a&1a&2a\\
b&2b&1b
\end{array}
\\
(a,b)(1,2)
\exend
\]
Here 
$A \cross C$ is the set of locations in a matrix with rows indexed by $C$
and columns indexed by $A$.
An entry $1a$ represents $(1,a) \in B \cross C$, etc.
The $(1,2)(a,b)$ underneath indicates a symmetry of $f$,
obtained by taking
$\alpha$ to be the identity, $\beta=(1,2)$, and $\gamma=(a,b)$.
Performing these substitutions yields
\[
\fabc=
\begin{array}{l|ll}
&x&y\\
\hline
b&2b&1b\\ 
a&1a&2a
\end{array}
\]
This is just a different representation of $f$,
as we see by swapping the rows,
so $\fabc=f$.
But we can't have $\hab=h$,
because $\alpha$ is the identity while $\beta$ is not,
so this $f$ has no $S(C)$-equivariant quotient,
hence $S(C)$ is not finitely cancelling.

We can simplify the display of this example as follows:
\[
\begin{array}{l|ll}
a&1&2\\
b&2&1
\end{array}
\]
\[
(a,b)(1,2)
\]
We don't need column labels as these aren't being permuted;
leaving out the labels from $C$ in the table entries
indicates this is a parallel bijection.

The example extends in an obvious way to show that $S(C)$ is not
finitely cancelling if $|C|>1$.
For example, take $C=\{a,b,c\}$, and
\[
\exstart
\begin{array}{l|lll}
a&1&2&3\\
b&2&3&1\\
c&3&1&2
\end{array}
\\
(a,b,c)(1,2,3)
\exend
\]

These examples come from the regular representation of a cyclic group.
A similar construction works for any finite group $G$.
(Cf. \ref{regrep} below.)
While we don't need it for what is to follow, we pause to illustrate
the construction in the case of the
noncyclic group $C_2 \cross C_2$,
whose regular representation is the Klein 4-group
$\{(a)(b)(c)(d),(a,b)(c,d),(a,c)(b,d),(a,d)(b,c)\}$:
\[
\exstart
\begin{array}{l|llll}
a&1&2&3&4\\
b&2&1&4&3\\
c&3&4&1&2\\
d&4&3&2&1
\end{array}
\\
(a,b)(c,d)(1,2)(3,4)
\\
(a,c)(b,d)(1,3)(2,4)
\exend
\]
This bijection is more symmetrical than we need to show this $\Gamma$ is
not cancelling,
because $\Gamma$ has a subgroup the two element subgroup generated by
$(a,b)(c,d)$,
and to show this is noncancelling we can just duplicate our first example
above:
\[
\exstart
\begin{array}{l|ll}
a&1&2\\
b&2&1\\
c&1&2\\
d&2&1
\end{array}
\\
(a,b)(c,d)(1,2)
\exend
\]

By now it is clear how to a handle any nontrivial permutation
all of whose cycles have the same length.
Such permutations are called \emph{semiregular}.
A permutation group is semiregular just if every non-trivial
element is semiregular.
(Such groups are also called `fixed point free', but this invites
confusion with groups with no globally fixed point.)
To sum up:

\begin{prop}[Feldman-Propp] \label{bad}
No permutation group
containing a semiregular subgroup
is finitely cancelling.
\end{prop}

Going further, Feldman and Propp give a beautiful algebraic proof of
the following:

\begin{theorem}
[Feldman-Propp]
\label{finitejustif}
A permutation group is finitely cancelling just if it has a globally
fixed point.
\end{theorem}

For further discussion, see \ref{finitecase} below.
For now, we're set:
We already have the tools to dispose of the infinite case.

\section{Not fully cancelling} \label{generalcase}

When $A$ and hence $B$ may be infinite,
known division methods depend on fixing an ordering for $C$.
This raises the suspicion that no nontrivial permutation group
can be fully cancelling.

\begin{theorem} \label{infinite}
A permutation group is fully cancelling
just if it is trivial.
\end{theorem}

In other words, if we demand complete equivariance for
$A$ and $B$, we can't demand any equivariance at all for $C$.

The proof will proceed via a string of examples.

We begin by slightly varying the construction used above in the finite case,
substituting non-parallel bijections.

\begin{itemize}
\item
$(a,b)$
\[
\exstart
\begin{array}{l|ll}
a&Ka&Kb\\
b&Qb&Qa
\end{array}
\\
(a,b)(K,Q)
\exend
\]

\item
$(a,b,c)$
\[
\exstart
\begin{array}{l|lll}
a&Ka&Kb&Kc\\
b&Qb&Qc&Qa\\
c&Jc&Ja&Jb
\end{array}
\\
(a,b,c)(K,Q,J)
\exend
\]

\item
$(a,b)(c,d)$ (not the simplest example; better for generalization)
\[
\exstart
\begin{array}{l|llll}
a&Ka&Kb&Kc&Kd\\
b&Qb&Qa&Qd&Qc\\
c&Ja&Jb&Jc&Jd\\
d&Xb&Xa&Xd&Xc
\end{array}
\\
(a,b)(c,d)(K,Q)(J,X)
\exend
\]

\end{itemize}

Now we jazz up these examples to include fixed points for the action
on $C$,
which we can't do in the finite case.

\begin{itemize}
\item
$(a,b)(c)$
\[
\exstart
\begin{array}{l|lllllll}
a&Ka&Kb&Kc&1a&2a&3a&\ldots\\
b&Qb&Qa&Qc&1b&2b&3b&\ldots\\
c&1c&2c&3c&4c&5c&6c&\ldots
\end{array}
\\
(a,b)(K,Q)
\exend
\]

\item
$(a,b,c)(d)$
\[
\exstart
\begin{array}{l|lllllllll}
a&Ka&Kb&Kc&Kd&1a&2a&3a&4a&\ldots\\
b&Qb&Qc&Qa&Qd&1b&2b&3b&4b&\ldots\\
c&Jc&Ja&Jb&Jd&1c&2c&3c&4c&\ldots\\
d&1d&2d&3d&4d&5d&6d&7d&8d&\ldots
\end{array}
\\
(a,b,c)(K,Q,J)
\exend
\]

\item
$(a,b,c)(d)(e)$
\[
\exstart
\begin{array}{l|llllllllll}
a&Ka&Kb&Kc&Kd&Ke&1a&2a&3a&4a&\ldots\\
b&Qb&Qc&Qa&Qd&Qe&1b&2b&3b&4b&\ldots\\
c&Jc&Ja&Jb&Jd&Je&1c&2c&3c&4c&\ldots\\
d&1d&2d&3d&4d&5d&6d&7d&8d&9d&\ldots\\
e&1e&2e&3e&4e&5e&6e&7e&8e&9e&\ldots
\end{array}
\\
(a,b,c)(K,Q,J)
\exend
\]

\item
$(a,b)(c,d)(e)$
\[
\exstart
\begin{array}{l|lllllllllll}
a&Ka&Kb&Kc&Kd&Ke&1a&2a&3a&4a&\ldots\\
b&Qb&Qa&Qd&Qc&Qe&1b&2b&3b&4b&\ldots\\
c&Ja&Jb&Jc&Jd&Je&1c&2c&3c&4c&\ldots\\
d&Xb&Xa&Xd&Xc&Xe&1d&2d&3d&4d&\ldots\\
e&1e&2e&3e&4e&5e&6e&7e&8e&9e&\ldots
\end{array}
\\
(a,b)(c,d)(K,Q)(J,X)
\exend
\]
\end{itemize}

These examples illustrate the method to prove that we can never require
any kind of equivariance for $C$.
The reason is that any nontrivial $\Gamma$ will contain some element that
is a product of one or more disjoint non-trivial cycles of the same length,
together with some fixed points.

\section{More about the regular representation} \label{regrep}

For future reference,
let's look more closely at the construction that we've been using,
based on the regular representation.

Fix a finite group $G$.
Take $A=B=C=G$,
\[
f = \{((x,y),(xy,y))\}
.
\]
(The unbound variables $x$ and $y$ are understood to range over $G$.)
First we observe that any quotient $h$ that is even $S_0(C)$-equivariant
will need to agree with one of the `rows' $f\row{c}$ of $f$.
To see this, fix $g \in G$ and set
\[
\alpha  = \beta = \{(x,gx)\}
.
\]
(The unbound variable $x$ is understood to range over $G$; you get the idea.)
Now
\[
f_{\alpha,\beta,\id}
= 
\{((gx,y),(gxy,y))\}
=
\{((x',y),(x'y,y))\}
= f
,
\]
so
\[
h(g) = h_{\alpha,\beta}(g) = \alpha^{-1} \then h \then \beta(g) 
= gh(g^{-1}g)
= gh(1)
= f\row{h(1)}(g)
.
\]
Since this holds for every $g \in G$,
\[
h = f\row{h(1)}
.
\]

Any row of $f$ will do as an $S_0(C)$-equivariant quotient, but
we can't have equivariance for any non-trivial element of $G$ acting on the
right.
Indeed, for any $g \in G$, we can take
\[
\beta = \gamma = 
\{(x,xg)\}
,
\]
\[
f_{\id,\beta,\gamma} = 
\{((x,yg),(xyg,yg))\}
= 
\{((x,g'),(xg',g'))\}
= f
.
\]
So we must have
\[
h = h_{\id,\beta} = h \then \beta
,
\]
that is,
\[
h(x) = h(x)g
,
\]
but this is impossible if $g$ is not the identity.

\section{Unfinished business}

\subsection{Back to the finite case} \label{finitecase}

Having determined exactly which groups $\Gamma \subset S_C$ are
fully cancelling,
we naturally turn our attention back to the finite case.
We've quoted Feldman and Propp's result
(Theorem \ref{finitejustif})
that $\Gamma$ is finitely cancelling just if it has a globally fixed point.
We've see that this condition is sufficient,
and shown that
if $\Gamma$
contains a fixed-point free subgroup
it is not finitely cancelling.
What about intermediate cases,
like the cyclic group generated by $(a,b,c)(d,e)$,
i.e. the group generated by $(a,b,c)$ and $(d,e)$,
where there is no fixed-point free subgroup?
Or the Klein-like 4-group
\[
\{\id,(a,b)(c,d),(a,b)(e,f),(c,d)(e,f)\}
,
\]
where there are no fixed-point free elements at all?
Feldman and Propp's beautiful algebraic proof does not immediately
provide counterexamples, though it gives a method to produce them.
They ask
\cite[Problem 4]{fp}
for more direct combinatorial arguments.

Let's at least dispose of $(a,b,c)(d,e)$:
\newcommand{\br}[1]{\bar #1}
\[
\exstart
\begin{array}{l|llllllllllll}
&\br0\br0&\br0\br1&\br0\br2&\br1\br0&\br1\br1&\br1\br2&00&01&02&10&11&12\\
\hline
a&\br00&\br01&\br02&\br10&\br11&\br12&0\br0&0\br1&0\br2&1\br0&1\br1&1\br2\\
b&\br01&\br02&\br00&\br11&\br12&\br10&0\br2&0\br0&0\br1&1\br2&1\br0&1\br1\\
c&\br02&\br00&\br01&\br12&\br10&\br11&0\br1&0\br2&0\br0&1\br1&1\br2&1\br0\\
d&0\br0&0\br1&0\br2&1\br0&1\br1&1\br2&\br00&\br01&\br02&\br10&\br11&\br12\\
e&1\br0&1\br1&1\br2&0\br0&0\br1&0\br2&\br10&\br11&\br12&\br00&\br01&\br02
\end{array}
\\
(a,b,c)
(\br00,\br01,\br02)(\br10,\br11,\br12)
(00,01,02),(10,11,12)\\
(d,e)
(0\br0,1\br0)(0\br1,1\br1)(0\br2,1\br2)
(00,10)(01,11)(02,12)\\
(a,b,c)(d,e)
(\br00,\br01,\br02)(\br10,\br11,\br12)
(0\br0,1\br0)(0\br1,1\br1)(0\br2,1\br2)
(00,11,02,10,01,12)
\exend
\]

This arises as follows.
Start with bijections
\[
p:X_1 \cross \{a,b,c\} \to X_2 \cross \{a,b,c\}
;\;
q:Y_1 \cross \{d,e\} \to Y_2 \cross \{d,e\}
,
\]
\[
\exstart
p =
\begin{array}{l|lll}
&\br0&\br1&\br2\\
\hline
a&0&1&2\\
b&1&2&0\\
c&2&0&1
\end{array}
\\
(\br0,\br1,\br2)(0,1,2)
\exend
\]
\[
\exstart
q =
\begin{array}{l|ll}
&\br0&\br1\\
\hline
d&0&1\\
e&1&0
\end{array}
\\
(\br0,\br1)(0,1)
\exend
.
\]
The inverses
\[
p^{-1}:X_2 \cross \{a,b,c\} \to X_1 \cross \{a,b,c\}
;\;
q^{-1}:Y_2 \cross \{d,e\} \to Y_1 \cross \{d,e\}
\]
are
\[
\exstart
p^{-1} =
\begin{array}{l|lll}
&0&1&2\\
\hline
a&\br0&\br1&\br2\\
b&\br2&\br0&\br1\\
c&\br1&\br2&\br0
\end{array},
\exend
\]
\[
\exstart
q^{-1} =
\begin{array}{l|ll}
&0&1\\
\hline
d&\br0&\br1\\
e&\br1&\br0
\end{array}
\exend
.
\]

Take the disjoint unions $X = X_1 \union X_2$ and $Y= Y_1 \union Y_2$
and augment $p$ and $q$ to involutions
\[
P = p \union p^{-1} \in S(X \cross \{a,b,c\});\;
Q = q \union q^{-1} \in S(Y \cross \{d,e\})
.
\]
Take products with the identity and combine to get an involution
\[
F=
P \cross id_Y
\union 
Q \cross id_X
\in
S(X \cross Y \cross \{a,b,c,d,e\})
.
\]
Let
\[
A = X_1 \cross Y_1 \union X_2 \cross Y_2
\]
and
\[
B = X_2 \cross Y_1 \union X_1 \cross Y_2
.
\]
Separate the involution $F$ into pieces
\[
F = f \union f^{-1}
,
\]
\[
f : A \cross \{a,b,c,d,e\} \to B \cross \{a,b,c,d,e\}
.
\]

This checkered Cartesian product construction
can be extended to cover any permutation without fixed points.
Any transitive permutation group contains such an element,
because the average number of fixed points is $1$,
and the identity has more.
So no transitive permutation group is finitely cancelling.

This construction
also takes care of our Klein-like 4-group.
In fact, it should handle
any subdirect product of nontrivial
cyclic permutation groups
(cf. Hall \cite[p.\ 63]{hall:groups}).
Now (asks Shikhin Sethi),
what about the 6-element group
\[
\Gamma = \{\id,(a,b,c), (b,c,a), (a,b)(d,e), (a,c)(d,e), (b,c)(d,e)\}
?
\]

\subsection{Deducing an ordering from a division method} \label{prob:reading}

A \emph{division method for $C$}
associates to any bijection
\[
f:A \cross C \to B \cross C
\]
a quotient bijection $Q(f)$ with the property that for any bijections
\[
\alpha:A \to A',\;\beta:B \to B'
,
\]
for the transformed division problem
\[
f_{\alpha,\beta}
=
(\alpha^{-1} \cross \id_C) \then f \then (\beta \cross \id_C):
A' \cross C \to B' \cross C
\]
the quotient
\[
Q(f_{\alpha,\beta}): A' \to B'
\]
satisfies
\[
Q(f_{\alpha,\beta})
=
Q(f)_{\alpha,\beta} = \alpha^{-1} \then Q(f) \then \beta
.
\]

A division method produces $S_0(C)$-equivariant quotients,
as we see by restricting $(\alpha,\beta)$ to $S(A) \cross S(B)$,
but more is required.
The method must not only respect symmetries of a particular problem,
it must give the same answer when presented with the same problem
in a different guise.
To see the distinction, consider that for an $f$ with no symmetries,
any bijection $h:A \to B$ is an $S_0(C)$-equivariant quotient,
and if a division method were required merely to respect the symmetries of $f$,
it could return a bijection depending on stupid properties of the set $A$,
like whether it consists entirely of natural numbers.

Once again we distinguish between full and finite division methods.
The method of Feldman and Propp is equivariant,
and yields finite division methods (one for each choice of basepoint in $C$).
In the infinite case we get division methods that depend on fixing an
ordering of $C$, and this dependence on the ordering
seems to be unavoidable.

\begin{problem}
Can we equivariantly associate a total ordering of $C$
to any full division method for $C$?
\end{problem}

In the finite case, we ask:

\begin{problem}
Can we equivariantly associate a single point in $C$
to any finite division method for $C$?
\end{problem}

The equivariance we're asking for here
means that we can't make arbitrary choices that favor one ordering
or point of $C$ over another.
Rather than fuss over the definition, let's consider the
particular case of division by three.

First, a general observation:
If $Q(f) = f\row{c}$ then
$Q(f_{\alpha,\beta}) = f_{\alpha,\beta}\row{c}$.
Indeed, for any $f$ we have
\[
(f\row{c})_{\alpha,\beta}
=
f_{\alpha,\beta}\row{c}
,
\]
so if $Q(f) = f\row{c}$,
\[
Q(f_{\alpha,\beta})
= Q(f)_{\alpha,\beta}
= (f\row{c})_{\alpha,\beta}
= f_{\alpha,\beta}\row{c}
.
\]

Now take $C=\{a,b,c\}$.
Consider the six bijections of the form
\[
\exstart
f[x,y,z] =
\begin{array}{l|lll}
&\br0&\br1&\br2\\
\hline
x&0&1&2\\
y&1&2&0\\
z&2&0&1
\end{array}
\\
(\br0,\br1,\br2)(0,1,2)
\exend
,
\]
where we propose to plug in for $x,y,z$
each of the six arrangements of $a,b,c$.
These six problems are really one and the same problem in six
different guises, because
\[
f[x,y,z]_{\id,(0,1,2)}
=
\begin{array}{l|lll}
&\br0&\br1&\br2\\
\hline
x&1&2&0\\
y&2&0&1\\
z&0&1&2\\
\end{array}
=
\begin{array}{l|lll}
&\br0&\br1&\br2\\
\hline
z&0&1&2\\
x&1&2&0\\
y&2&0&1\\
\end{array}
= f[z,x,y]
,
\]
and
\begin{eqnarray*}
f[x,y,z]_{(\br1,\br2),(1,2)}
&=&
\begin{array}{l|lll}
&\br0&\br2&\br1\\
\hline
x&0&2&1\\
y&2&1&0\\
z&1&0&2
\end{array}
\\&=&
\begin{array}{l|lll}
&\br0&\br1&\br2\\
\hline
x&0&1&2\\
y&2&0&1\\
z&1&2&0\\
\end{array}
\\&=&
\begin{array}{l|lll}
&\br0&\br1&\br2\\
\hline
x&0&1&2\\
z&1&2&0\\
y&2&0&1
\end{array}
\\&=&
f[x,z,y]
.
\end{eqnarray*}
A division method must produce a quotient respecting the symmetry
\[
f[x,y,z]_{(\br0,\br1,\br2),(0,1,2)} = f[x,y,z]
,
\]
so it must conjugate the cycle $(\br0,\br1,\br2)$ to the cycle $(0,1,2)$.
There are three ways to do this, corresponding to the three
rows $x,y,z$ in the table,
so (as observed in section \ref{regrep}) the quotient bijection
$Q(f[x,y,z])$
distinguishes one of the three elements of $C$, which we call $\bp[x,y,z]$:
\[
Q(f[x,y,z]) = f[x,y,z]\row{\bp[x,y,z]}
.
\]
By the general result above, these six basepoints $\bp[x,y,z]$ all coincide.
So we can distinguish a basepoint in $\bp \in C$
without making any arbitrary choices
of how to order the elements of $C$.
This is the kind of equivariance we're looking for.

For a finite division method, that's as far as we can go.
In the infinite case,
say that our distinguished basepoint $\bp$ is $c$.
We continue by presenting the two problems $f[a,b],f[b,a]$, where
\[
\exstart
f[x,y]=
\begin{array}{l|lllllll}
&\br1&\br2&\br3&\br4&\br5&\br6&\ldots\\
\hline
x&Kx&Ky&Kc&1x&2x&3x&\ldots\\
y&Qy&Qx&Qc&1y&2y&3y&\ldots\\
c&1c&2c&3c&4c&5c&6c&\ldots
\end{array}
\exend
.
\]
The bijection $f[x,y]$
in effect associates $K$ with $x$ and $Q$ with $y$;
depending on where $K$ and $Q$ wind up
under the quotient bijection $Q(f[x,y])$
(or rather, its inverse),
we can pick $K$ over $Q$, hence $x$ over $y$.
Our preference of $a$ over $b$ will be the same whether we use
$f[a,b]$ or $f[b,a]$,
because these are really the same problem:
\[
f[x,y]_{\id,(K,Q)} = f[y,x]
.
\]

Now, what about division by four? Or five?

\subsection{Parallelizing a bijection} \label{prob:parallel}

We've already observed that while $\Gamma$ is finitely cancelling
just if every parallel bijection has an equivariant quotient,
if $\Gamma$ is not finitely cancelling there could be special
bijections $f$ which have a $\Gamma$-equivariant quotient,
while their Feldman-Propp parallelizations $\bar{f}$ do not.

\begin{problem}
If a finite bijection $f: A \cross C \to B \cross C$
has a $\Gamma$-equivariant quotient, must the parallelization
$\bar{f}$ also have a $\Gamma$-equivariant quotient?
\end{problem}

We haven't thought very hard about this one.

\subsection{Special cases}

There are plenty of other questions we could ask,
say concerning restrictions that will guarantee
that $S(C)$-equivariant division is possible.
For example, we might fix $n,k$ and
ask whether $S(C)$-equivariant division is always
possible
when $|A|=|B|=n$ and $|C|=k$.
It is easy to see that in this case we must have $\gcd(k,n!)=1$,
i.e.\ $k$ must have no prime factor $\leq n$.
This condition is sufficient for $n=1,2,3$ and maybe $4$;
the proofs get more involved as $n$ increases.
On the other hand, an example
(thanks to John Voight)
shows that division is not
always possible when $n=8$ and $k=11$.

\section*{Thanks}
Thanks to David Feldman and Shikhin Sethi for crucial advice.

\bibliography{equi}
\bibliographystyle{hplain}

\end{document}